\documentclass{amsart}
\usepackage{amssymb}
\usepackage[centertags]{amsmath}
\usepackage{amssymb}
\usepackage{xr}
\usepackage{amscd}
\usepackage{amsthm}
\usepackage{amsxtra}

\newcommand{\envelope}[1]{\ensuremath{\widetilde{#1}} }
\newcommand{\re}{\ensuremath{\text{Re }}}
\newcommand{\im}{\ensuremath{\text{Im }}}

\theoremstyle{definition}
\newtheorem{theorem}{Theorem}[section]
\newtheorem{lemma}[theorem]{Lemma}
\newtheorem{cla}[theorem]{Claim}
\newtheorem{problem}[theorem]{Problem}
\newtheorem{definition}[theorem]{Definition}
\newtheorem{remark}[theorem]{Remark}
\newtheorem{prop}[theorem]{Proposition}

\begin{document}
\title{A schlichtness theorem for envelopes of holomorphy}

\author{Daniel Jupiter}
\address{Department of Mathematics, Texas A \& M University,
College Station TX, 77843-3368, U. S. A.}
\email{jupiter@math.tamu.edu}

\begin{abstract}
Let $\Omega$ be a domain in $\mathbb{C}^2$. We prove the following
theorem. If the envelope of holomorphy of $\Omega$ is schlicht
over $\Omega$, then the envelope is in fact schlicht. We provide
examples showing that the conclusion of the theorem does not hold
in $\mathbb{C}^n$, $n>2$. Additionally, we show that the theorem
cannot be generalized to provide information about domains in
$\mathbb{C}^2$ whose envelopes are multiply sheeted.
\end{abstract}

\maketitle

\section{Introduction}
\label{section:intro} A distinct difference between complex
analysis in one variable and in several variables is the existence
in higher dimension of domains which are not Stein. We are led to
consider the envelope of holomorphy of such a domain: the smallest
Stein Riemann domain over $\mathbb{C}^n$ to which all functions
holomorphic on our domain extend holomorphically. (See, for
example \cite[Chapter H]{gunning:scv}.)

In the case of domains in $\mathbb{C}^n$ with relatively simple
structure, the envelope of holomorphy can, at least in principle,
be calculated explicitly. Examples of such domains are Reinhardt,
Hartogs, tubular and contractible domains. These classes of
domains enjoy the distinct advantage that their envelopes are
schlicht; that is, their envelopes are domains in $\mathbb{C}^n$.

In general, however, the envelope of a domain need not be
schlicht. It may be a finitely or infinitely sheeted Riemann
domain spread over $\mathbb{C}^n$. It is not obvious, given a
description of a domain, whether its envelope is schlicht or
multiply sheeted. Nor is it entirely clear how the geometry of a
given domain relates to the geometry of its envelope. In fact, we
have the following open question.

\begin{problem}
Does there exist a bounded domain, $\Omega$, in $\mathbb{C}^n$,
$n\geq 2$, such that $\Omega$ has smooth boundary, and the
envelope of $\Omega$ is infinitely sheeted?
\end{problem}

There are classical examples which show that there exist bounded
domains whose envelopes are infinitely sheeted. Unfortunately,
these domains have boundaries which are highly nonsmooth. (Such an
example can be constructed, e.g., by nesting infinitely many
Hartogs figures one inside the next, and carefully connecting
them.)\\

An obvious problem suggested by the above discussion is that of
determining why and when the envelope is schlicht or multiply
sheeted. Can we, just by looking at a domain, determine whether
its envelope is schlicht? Is there a way of controlling or
understanding the number of sheets in the envelope of a domain?
While we are unable to fully answer these questions, we present a
theorem which indicates that in certain situations there are
criteria which allow us to better understand whether envelopes are
schlicht or not.

Specifically, we prove the following.
\begin{theorem}\label{theorem:main theorem} Let $\Omega$ be a domain
in $\mathbb{C}^2$, and let $(\widetilde{\Omega},\,\pi)$ be its
envelope of holomorphy. If $\pi$ is injective on
$\pi^{-1}(\Omega)$, then $\pi$ is injective on
$\widetilde{\Omega}$. In other words, if $\widetilde{\Omega}$ is
schlicht over $\Omega$, then $\widetilde{\Omega}$ is schlicht.
\end{theorem}

The theorem indicates that for domains in $\mathbb{C}^2$ there may
be a relationship between the number of sheets in the envelope
lying above the domain, and the number of sheets in the envelope.
As we shall see in Section \ref{section:counterexamples}, there
appears to be a relationship only in the special case of domains
whose envelope is schlicht over the domain. We shall also see in
Section \ref{section:counterexamples} that there is no such
relationship in higher dimensions. Further examples in this vein
can be found in \cite{jupiter:thesis}.\\

We note that Chirka and Stout proved a weaker result in
\cite{chirka-stout:schlicht}. Specifically, they proved Lemma
\ref{lemma:lemma 1}. Our methods are different, and better suited
to proving and understanding the results in which we
are interested.\\

\section{A Schlichtness Theorem}
\label{section:theorem} We prove Theorem \ref{theorem:main
theorem}. Our proof proceeds as follows. Using the natural
embedding of $\Omega$ into $\envelope{\Omega}$, we are considering
$\Omega$ as a subset of $\envelope{\Omega}$ rather than as a
domain in $\mathbb{C}^2$. Let $\Omega_1$ be the set of points in
$\widetilde{\Omega}$ which can be reached from $\Omega$ by pushing
discs. Similarly, let $\Omega_n$ be the set of points in
$\widetilde{\Omega}$ which can be reached from $\Omega_{n-1}$ by
pushing discs. We recall that if two one-dimensional analytic
varieties in $\mathbb{C}^2$ intersect nontrivially, then so do
slight perturbations of these varieties. We proceed inductively,
showing that $\pi$ is injective on $\Omega_n$ for each $n$. In
other words, we show that for each $n$, $\Omega_n$ can be
identified with a domain in $\mathbb{C}^2$. To do so we push discs
from $\Omega_{n-1}$, keeping track of their intersections. Using
the above fact about analytic varieties, we show that if
$\Omega_n$ is not schlicht, then it is not schlicht over
$\Omega_{k}$, $k\leq n-1$. In particular $\Omega_n$ is not
schlicht over $\Omega$, contradicting the hypothesis that
$\envelope{\Omega}$ is schlicht over $\Omega$.

We make several definitions which we shall need in the proof of
the theorem.

\begin{definition}[Pushing Discs]\label{definition:pushing discs}
We say that a point, $p\in\envelope{\Omega}$, can be reached from
$\Omega$ by pushing discs if there is a neighbourhood $U$ of $p$
in $\envelope{\Omega}$ such that $\pi|_{U}$ is a biholomorphism,
and such that the following holds: there is a biholomorphism, $F$,
of $\Delta^2(0,\,1)$ into $U$ such that
\begin{enumerate}
\item $p\in F(\Delta^2(0,\,1))$, and
\item $F(H)\subset\subset U\cap\Omega$,
\end{enumerate}
where $H$ is the Hartogs figure,
\[H=\biggl(\Delta(0,\,1)\times
\biggl\{\frac{1}{2}<|w|<1\biggr\}\biggr)
\cup\biggl(\Delta\biggl(0,\,\frac{1}{2}\biggr)\times\Delta(0,\,1)\biggr).\]
\end{definition}

\begin{remark} Several observations and remarks about disc pushing should be made.\\
\begin{enumerate}
\item We insist that our neighbourhoods, $U$, are biholomorphic to
balls in $\mathbb{C}^2$, in particular $\pi(U)$ is a ball. \item
The extension of a holomorphic function from $F(H)$ to
$F(\Delta^2(0,\,1))$ is single valued. \item Given a point,
$p\in\envelope{\Omega}$, which can be reached from $\Omega$ by
pushing discs as above, we can assume (by rotation and scaling in
the first coordinate, if necessary) that $p$ is in $F(C)$, where
\[C=[0,\,1]\times\Delta(0,\,1).\]
In other words, we have a continuous family of holomorphic discs
in $\envelope{\Omega}$ such that the ``bottom'' disc and the
boundaries of the discs lie in $\Omega$, while the ``top'' disc is
not contained in $\Omega$..
\end{enumerate}
\end{remark}

\begin{definition}[$\Omega_n$]\label{definition:omega n}
We inductively define the sets $\Omega_{n}$. Let
\[\Omega_0=\Omega.\]
Let
\begin{align*}
\Omega_{n+1}= \{&\text{ all points in } \envelope{\Omega} \text{
which can be reached from }\\ & \text{ }\Omega_n \text{ by pushing
discs }\}.
\end{align*}
\end{definition}

\begin{remark} We make several observations about $\Omega_n$.

\begin{enumerate}
\item It is clear from the construction of $\Omega_n$ that these
sets are open subsets of $\envelope{\Omega}$. \item We have
defined $\Omega_n$ as a subset of $\envelope{\Omega}$, and when
building $\Omega_n$ we push discs within $\envelope{\Omega}$.
However, if $\Omega_{n-1}$ is in fact a domain in $\mathbb{C}^2$,
we can also view disc pushing as follows. We have
\[\pi(F(\Delta^2(0,\,1)))\subset
\mathbb{C}^2, \text{ where }
\pi(F(H))\subset\pi(\Omega_{n-1})=\Omega_{n-1}.\] In other words
we push discs from $\Omega_{n-1}$, considering $\Omega_{n-1}$ as a
subset of $\mathbb{C}^2$. We then lift these discs to
$\Omega_{n}\subset\envelope{\Omega}$. By the identity principle
for liftings (\cite[Proposition 1.1.5]{jarnicki-pflug:extension})
the liftings are unique.
\end{enumerate}
\end{remark}

We next observe that we can use $\Omega_n$ to recover $\envelope{\Omega}$.

\begin{prop}
Let $\Omega'=\cup_{n=0}^{\infty}\Omega_n$. Then
$\Omega'=\widetilde{\Omega}$.
\end{prop}

\begin{proof} Assume not, and let $p\in\envelope{\Omega}$ be a point
in $\partial\Omega'$. We claim that there is a  neighbourhood $U$
of $p$ in $\envelope{\Omega}$ such that $U\cap\Omega'$ is
pseudoconvex.

In fact, let $U$ be a neighbourhood of $p$ such  that $\pi|_U$ is
a biholomorphism, and further choose $U$ so that $\pi(U)$ is a
ball. If $U\cap\Omega'$ is not pseudoconvex, then we can find a
biholomorphism, $F$, of $\Delta^2(0,\,1)$ into $U$ such that
\begin{enumerate}
\item $F(H)\subset\subset U\cap\Omega'$, and
\item $F(\Delta^2(0,\,1))$ contains points which are not in $\Omega'$.
\end{enumerate}
This follows from the fact that $\pi|_U$ is a biholomorphism, so
that we can consider $U\cap\Omega'$ as a domain in $\mathbb{C}^2$,
and from the notion of II-pseudoconvexity as described in
\cite[Chapter II.2]{pflug:pc}.

Now $F(H)\subset\subset\Omega'$, so in fact
$F(H)\subset\subset\Omega_n$, for some $n$. This implies that
$\Omega_{n+1}$ contains points which are not in $\Omega'$. This is
a contradiction.

We conclude that every point in $\partial\Omega'$ has a
neighbourhood, $U$, such that $U\cap\Omega'$ is pseudoconvex. By
the equivalence of local and global pseudoconvexity of unbranched
Riemann domains over $\mathbb{C}^n$ (see e.g. \cite[Corollary
2.2.16]{jarnicki-pflug:extension}) we see that $\Omega'$ is
pseudoconvex. However, as $\Omega'$ is a subset of
$\envelope{\Omega}$, every holomorphic function on $\Omega$
extends in a single valued fashion to $\Omega'$. We conclude that
$\Omega'=\envelope{\Omega}$.\qed
\end{proof}

To prove Theorem \ref{theorem:main theorem} we require two lemmas.

\begin{lemma}\label{lemma:lemma 1}
Assume that $\Omega_n$ is schlicht. Assume on the other hand that
$\Omega_{n+1}$ is not schlicht. Then $\Omega_{n+1}$ is not
schlicht over $\Omega_n$.

Precisely, assume that there are points, $p_1\neq p_2$ in
$\Omega_{n+1}$ such that $\pi(p_1)=\pi(p_2)$. Then in fact there
are points $q_1\neq q_2$ in $\Omega_{n+1}$ such that
$q=\pi(q_1)=\pi(q_2)\in\Omega_n$.
\end{lemma}

This lemma says, in particular, that if pushing discs from
$\Omega$ once does not create any sheets over $\Omega$, then in
fact it does not create any sheets over $\Omega_1$. As noted in
the Introduction, Chirka and Stout have also proved this lemma
\cite{chirka-stout:schlicht}.

\begin{lemma}\label{lemma:lemma 2}
Assume that $\Omega_n$ is schlicht. If $\Omega_{n+1}$ is not
schlicht over $\Omega_k$, $k\leq n$, then $\Omega_{n+1}$ is not
schlicht over $\Omega_{k-1}$.
\end{lemma}

Once we have these lemmas we prove the theorem as follows.

\begin{proof}[Proof of Theorem \ref{theorem:main theorem}]
We inductively show that $\Omega_n$ is schlicht. By assumption we
have that $\Omega_0$ is schlicht.

If $\Omega_1$ were not schlicht, then Lemma \ref{lemma:lemma 1}
would imply that it was not schlicht over $\Omega_0$. This in turn
means that $\envelope{\Omega}$ is not schlicht over $\Omega_0$.
This is a contradiction, and we conclude that $\Omega_1$ is
schlicht.

Applying Lemma \ref{lemma:lemma 1} to $\Omega_2$ reveals that if
$\Omega_2$ were not schlicht, then it would not be schlicht over
$\Omega_1$. Lemma \ref{lemma:lemma 2} now implies that $\Omega_2$
is not schlicht over $\Omega_0$. This again gives a contradiction,
and we conclude that $\Omega_2$ is schlicht.

We proceed similarly for $\Omega_n$, $n\geq 3$.\qed
\end{proof}

We now prove Lemmas \ref{lemma:lemma 1} and \ref{lemma:lemma 2}.

\begin{proof}[Proof of Lemma \ref{lemma:lemma 1}]
We begin by choosing two points, $p_1\neq p_2$ in $\Omega_{n+1}$
with $p=\pi(p_1)=\pi(p_2)$. By definition and our comments after
Definition \ref{definition:pushing discs}, for each point $p_j$
there is a neighbourhood, $U_j$, of $p_j$ in $\envelope{\Omega}$,
and a biholomorphism $F_j$ of $\Delta^2(0,\,1)$ into
$\envelope{\Omega}$ satisfying
\begin{enumerate}
\item $p_j\in F_j(C)$, and \item $F_j(H)\subset\subset
U_j\cap\Omega_n$,
\end{enumerate}

recalling that $C$ is defined as
\[C=[0,\,1]\times\Delta(0,\,1).\] We now define a subset $c$ of $C$ as
\[c=(\{0\}\times\Delta(0,\,1))\cup
([0,\,1]\times\partial\Delta(0,\,1)).\]

We have, by assumption, that $\pi(F_1(C))\cap\pi(F_2(C))$ is
nonempty; in particular it contains $p$. We let $M$ be the
connected component of $\pi(F_1(C))\cap\pi(F_2(C))$ which contains
$p$. We note that the $\pi(F_j(C))$ can be assumed to be in
general position, so that their intersection is a two dimensional
manifold. (The intersection is a closed set, so strictly speaking
the intersection is a manifold with ``ends''.) By the stability of
intersection of one-dimensional analytic manifolds in
$\mathbb{C}^2$, we conclude that $M$ contains a point, $q$, which
lies in $\pi(F_1(c))$ or in $\pi(F_2(c))$. We then have that
$q=\pi(F_1(z_1))$ for some $z_1$ in $c$, or $q=\pi(F_2(z_2))$ for
some $z_2$ in $c$. In the first case $q_1=F_1(z_1)$ is in
$\Omega_n$, in the second case $q_2=F_2(z_2)$ is in $\Omega_n$.

We eliminate the possibility that both $q_1=F_1(z_1)$ and
$q_2=F_2(z_2)$ are points in $\Omega_n$. Otherwise, since
$\Omega_n$ is schlicht and these two points have the same
projection, we must have that $q_1=q_2$. Let $\gamma$ be a path in
$M$ connecting $q$ to $p$. By the identity principle for liftings,
we conclude that $p_1=p_2$. This is a contradiction.

We examine the case where $z_1\in c$, hence $F_1(z_1)$ is in
$\Omega_n$. Since $q_2=F_2(z_2)$ is not in $\Omega_n$, $q_1\neq
q_2$. But $\pi(q_1)=\pi(q_2)$. We now conclude that $\Omega_{n+1}$
is not schlicht over $\Omega_n$.

The case where $q$ is in $\pi(F_2(c))$ but not in $\pi(F_1(c))$
proceeds analogously.\qed
\end{proof}

The proof of Lemma \ref{lemma:lemma 2} is similar.

\begin{proof}[Proof of Lemma \ref{lemma:lemma 2}]
We begin by choosing two points, $p_1\neq p_2$, with $p_1$ in
$\Omega_{n+1}$ and $p_2$ in $\Omega_{k}$, and with
$p=\pi(p_1)=\pi(p_2)$. By definition, for each point $p_j$ there
is a neighbourhood, $U_j$, of $p_j$ in $\envelope{\Omega}$, and a
biholomorphism $F_j$ of $\Delta^2(0,\,1)$ into $\envelope{\Omega}$
satisfying
\begin{enumerate}
\item $p_1\in F_1(C)$, $F_1(H)\subset\subset U_1\cap\Omega_n$, and
\item $p_2\in F_2(C)$, $F_2(H)\subset\subset U_2\cap\Omega_{k-1}$.
\end{enumerate}

As in the proof of Lemma \ref{lemma:lemma 1}, we have that
$\pi(F_1(C))\cap\pi(F_2(C))$ is nonempty; in particular it
contains $p$. We let $M$ be the connected component of
$\pi(F_1(C))\cap\pi(F_2(C))$ which contains $p$. By the stability
of intersection of one-dimensional analytic manifolds in
$\mathbb{C}^2$, we conclude that $M$ contains a point, $q$, which
lies in $\pi(F_1(c))$ or in $\pi(F_2(c))$.

It is not possible that $F_1(z_1)$ is in $\Omega_n$. Indeed, since
$q_2=F_2(z_2)$ is in $\Omega_k\subset\Omega_n$, this would imply
that we have $F_1(z_1)=F_2(z_2)$ by the fact that $\Omega_n$ is
schlicht. As above, the path $\gamma\subset M$ connecting $q$ and
$p$ yields the contradiction that $p_1=p_2$.

We thus must have that $q_1=F_1(z_1)$ is in $\Omega_{n+1}$ but not
in $\Omega_n$, and that $q_2=F_2(z_2)\in F_2(c)$ and hence it is
in $\Omega_{k-1}$. This means that $\Omega_{n+1}$ is not schlicht
over $\Omega_{k-1}$. With this, the proof of the lemma is
complete.\qed
\end{proof}

\subsection{Remarks}\label{subsection:remarks}
Theorem \ref{theorem:main theorem} also holds for two
dimensional Riemann domains over $\mathbb{C}^2$. The proof proceeds exactly as above.\\

\section{Counterexamples} \label{section:counterexamples}
\subsection{Counterexamples in
$\mathbb{C}^2$}\label{subsection:examples in C2} We construct a
domain, $\Omega$, in $\mathbb{C}^2$ such that the envelope of
$\Omega$ has two sheets over $\Omega$, but three sheets over
$\mathbb{C}^2$. This shows that Theorem \ref{theorem:main theorem}
cannot be generalized to give information about domains in
$\mathbb{C}^2$ whose
envelopes are multiply sheeted.\\

The example is obtained as follows. We construct three families of
analytic discs whose intersection is a small set. For each family
we find a pseudoconvex domain close to the family, making sure
that certain boundary points of the new domain are in fact
strictly pseudoconvex. We call this new domain a fattening of the
family. We then build a frame for the family: a small
neighbourhood of the union of the bottom disc of the family and
the boundaries of the discs in the family. Pushing discs in this
frame gives us the fattened family with which we started. We join
the frames with well chosen paths and find pseudoconvex
neighbourhoods of the paths. The domain, $\Omega$, is the union of
the frames and the paths. The correct choice of paths ensures that
our domain has the desired properties.\\

We begin by constructing three families of discs,
\begin{alignat*}{2}
\Sigma_s:&\Delta(0,\,2)\rightarrow\mathbb{C}^2,&\quad
s\in(-1,\,\epsilon_{\Sigma}),\\
\Delta_t:&\Delta(0,\,1/2)
\rightarrow \mathbb{C}^2,&\quad t\in (-\epsilon_{\Delta},\,1),\\
\Gamma_l:&\Delta(0,\,r_{\Gamma})\rightarrow \mathbb{C}^2,&\quad
l\in(-\epsilon_{\Gamma},\,\epsilon_{\Gamma})
\end{alignat*}
where $\epsilon_{\Sigma}$, $\epsilon_{\Delta}$,
$\epsilon_{\Gamma}$ and $r_{\Gamma}$ are small positive real
numbers to be carefully chosen, and the maps are defined as
\begin{align*}
\Sigma_s(w) & =(s,\,w),\\ \Delta_t(z)&=(z,\,t),\\
\Gamma_l(\xi)&=(\xi+i\eta\xi^2,\,\xi-il),
\end{align*}
with $\eta$ a small positive real number to be carefully chosen.

We shall refer to the families as $\Sigma$, $\Delta$ and $\Gamma$.
By the bottom disc of the family $\Sigma$, $\Delta$ or $\Gamma$ we
mean, respectively, the discs $\Sigma_{-1}$, $\Delta_{1}$ or
$\Gamma_{\epsilon_{\Gamma}}$. For ease of notation we denote these
discs by $\Sigma_B$, $\Delta_B$ and $\Gamma_B$ respectively. By
the boundary of the family $\Sigma$, $\Delta$ or $\Gamma$ we mean,
respectively, the family of circles
\begin{align*}
\Sigma_s|_{\partial\Delta(0,\,2)},& \quad
s\in(-1,\,\epsilon_{\Sigma}),\\
\Delta_t|_{\partial\Delta(0,\,1/2)},&\quad
t\in(-\epsilon_{\Delta},\,1),\\
\Gamma_l|_{\partial\Delta(0,\,r_{\Gamma})}&\quad
l\in(-\epsilon_{\Gamma},\,\epsilon_{\Gamma}).
\end{align*}
We denote these boundaries by $\partial\Sigma$, $\partial\Delta$,
and $\partial\Gamma$, respectively. Finally,
$\Sigma'=\Sigma\cup\partial\Sigma$,
$\Delta'=\Delta\cup\partial\Delta$, and
$\Gamma'=\Gamma\cup\partial\Gamma$. Let $\pi_j$ be the projection
onto the $j$th coordinate of $\mathbb{C}^2$.\\

For each of the three following claims we must make suitable
choices of $\epsilon_{\Sigma}$, $\epsilon_{\Delta}$,
$\epsilon_{\Gamma}$, $r_{\Gamma}$ and $\eta$.

We would first like to show that the intersection of the three
families, $\Sigma\cap\Delta\cap\Gamma$, is a small set. We let
\begin{align*}P=\biggl\{(z_1,\,z_2)\in\mathbb{C}^2\,;\,
\im z_1&=\im z_2=0,\\
&\frac{-1}{2}\leq \re z_1 \leq
\epsilon_{\Sigma},\,-\epsilon_{\Delta}\leq \re z_2 \leq
1\biggr\}.\end{align*}
We see that
\[\Sigma'\cap\Delta'=P.\]

Let $G$ be the union of the graphs of the functions
\[f(y)=\pm\biggl(y^2-\frac{y}{\eta}\biggr)^{1/2},\,y\leq 0.\]
In other words
\[G=\{(x,\,y)\in\mathbb{R}^2\,;\,(f(y),\,y),\,y\leq 0\}.\]

A point in $\Gamma_l$ looks like $(\xi+i\eta\xi^2,\,\xi-il)$. Let
$\xi=x+iy$.

Choosing $r_{\Gamma}<1/\eta$ we see that for a point,
$\Gamma_l(\xi)=(\xi+i\eta\xi^2,\,\xi-il)$, in $\Gamma_l$ to be in
$P$ we must have that
\begin{enumerate}
\item $y\leq 0$,
\item $y=l$, and
\item $x=\pm\biggl(y^2-\frac{y}{\eta}\biggr)^{1/2}$.
\end{enumerate}
Thus $\Gamma_l(\xi)=\Gamma_y(\xi)$. If $l>0$ and
$\xi\in\Delta(0,\,r_{\Gamma})$ then $\Gamma_l(\xi)$ is not in $P$.
If $l=0$ and $\Gamma_l(\xi)$ is in $P$, then $\xi=0$. Finally, if
$l<0$ and $\Gamma_l(\xi)$ is in $P$, then
$\xi=\pm(l^2-l/\eta)^{1/2}+il$.

We have shown:
\begin{cla}\label{claim:intersections 1} The intersection of the
three families, $\Sigma\cap\Delta\cap\Gamma$, is a small set.
Specifically,
\[\Sigma\cap\Delta\cap\Gamma=P\cap\{\Gamma_y(\xi)\,;\,\xi=x+iy\in
G\cap\Delta(0,\,r_{\Gamma}),\, -\epsilon_{\Gamma}\leq y\leq 0\}.\]
\end{cla}

The following claim will allow us to choose the frames of our
families to be disjoint.

\begin{cla}\label{claim:intersections 2}
The boundaries of the families are pairwise disjoint. Similarly,
the bottom discs of the families are pairwise disjoint. The
boundary of each family is disjoint from the bottom disc of the
other families.
\end{cla}

Finally, we note that the boundaries and bottoms of our families
of discs are disjoint from the intersection of the three families.

\begin{cla}\label{claim:intersections 3}
The sets $\Sigma_B$, $\Delta_B$, $\Gamma_B$, $\partial\Sigma$,
$\partial\Delta$ and $\partial\Gamma$ are all disjoint from
$\Sigma\cap\Delta\cap\Gamma$.
\end{cla}

Given a family of discs such as $\Delta$, and any $\delta>0$, we
can find a pseudoconvex domain, $\Delta^{\Box}$, contained in a
$\delta$ neighbourhood of $\Delta$. We can then construct a
subdomain, $\Delta^{\sqcup}$, of $\Delta^{\Box}$ such that
\begin{enumerate}
\item $\Delta^{\sqcup}$ is contained within a $\delta$
neighbourhood of $\Delta_B\cup\partial\Delta$, \item
$\widetilde{\Delta}^{\sqcup}=\Delta^{\Box}$, and \item the
boundary of $\Delta^{\sqcup}$ contains smooth points which consist
of strictly pseudoconvex points.
\end{enumerate}
We call $\Delta^{\Box}$ a fattening of $\Delta$, and we call
$\Delta^{\sqcup}$ the frame of $\Delta$.

An analogous construction can be carried out for $\Gamma$ and
$\Sigma$. In fact, by Claim \ref{claim:intersections 2} we can
choose these frames to be pairwise disjoint. By Claim
\ref{claim:intersections 3} we can also choose them so that they
do not intersect
$\Sigma^{\Box}\cap\Delta^{\Box}\cap\Gamma^{\Box}$.\\

Let $\gamma_1$ and $\gamma_2$ be two disjoint smooth paths,
$\gamma_1:[0,\,1]\rightarrow\mathbb{C}^2$ and
$\gamma_2:[0,\,1]\rightarrow\mathbb{C}^2$, satisfying:
\begin{enumerate}
\item The path $\gamma_1$ connects a strictly pseudoconvex boundary
point of $\Delta^{\sqcup}$ and a strictly pseudoconvex boundary
point of $\Sigma^{\sqcup}$. Similarly, $\gamma_2$ connects a
strictly pseudoconvex boundary point of $\Sigma^{\sqcup}$ and a
strictly pseudoconvex boundary point of $\Gamma^{\sqcup}$.
\item The intersection of $\Sigma^{\Box}\cup\Delta^{\Box}\cup\Gamma^{\Box}$
and $\gamma_1\cup\gamma_2$ consists of the four endpoints of the
paths.
\item The curve $\gamma_1$ is transversal
to the boundaries of $\Delta^{\sqcup}$ and $\Sigma^{\sqcup}$ at
its endpoints. The analogous statement holds for $\gamma_2$.
\item Let $p$ be a point in the complement of $\pi_1(\Sigma^{\Box}\cup\Delta^{\Box}
\cup\Gamma^{\Box})$, and define the set $A$ as
\[A=\mathbb{C}^2\backslash(\{z_1=p\}\times\mathbb{C}).\]
The fundamental group of $A$ is nontrivial. We choose $\gamma_1$
in such a way that
\begin{enumerate}
\item $\Delta^{\sqcup}\cup\gamma_1\cup\Sigma^{\sqcup}$ contains no
nontrivial element of the fundamental group, and
\item $\Delta^{\Box}\cup\gamma_1\cup\Sigma^{\Box}$ does contain a
nontrivial element of the  fundamental group.
\end{enumerate}
We choose $\gamma_2$ so that analogous statements hold for
$\gamma_2$, $\Sigma$ and $\Gamma$.
\end{enumerate}

Item 4 implies that
$f(z_1,\,z_2)=(z_1-p)^{1/3}$ is holomorphic on a small neighbourhood of
\[\Sigma^{\sqcup}\cup\Delta^{\sqcup}\cup\Gamma^{\sqcup}\cup\gamma_1\cup\gamma_2,\]
but triple valued on a small neighbourhood of
\[\Sigma^{\Box}\cup\Delta^{\Box}\cup\Gamma^{\Box}\cup\gamma_1\cup\gamma_2.\]
Indeed, as we follow $\gamma_1$ from $\Delta^{\Box}$ to $\Sigma^{\Box}$ we
change branches of the cube root function. Similarly, as we follow
$\gamma_2$ from $\Sigma^{\Box}$ to $\Gamma^{\Box}$ we change
branches again.

By a theorem of Forn{\ae}ss and Stout
\cite{fornaess-stout:spreading-polydiscs} we can find a
neighbourhood, $\Gamma_1$, of $\gamma_1$ such that $\Gamma_1$ is
pseudoconvex, $\Gamma_1\cup\Delta^{\sqcup}$ is locally
pseudoconvex near the intersection of $\gamma_1$ and
$\Delta^{\sqcup}$, and $\Gamma_1\cup\Sigma^{\sqcup}$ is locally
pseudoconvex near the intersection of $\gamma_1$ and
$\Sigma^{\sqcup}$. Similarly, we find a neighbourhood, $\Gamma_2$,
of $\gamma_2$ with analogous properties. These neighbourhoods can be
made as small as we like.

We define $\Omega$ as
\[\Omega=\Sigma^{\sqcup}\cup\Delta^{\sqcup}\cup\Gamma^{\sqcup}\cup\Gamma_1\cup\Gamma_2.\]

As noted above, $f(z_1,\,z_2)=(z_1-p)^{1/3}$ is holomorphic on
$\Omega$, but triple valued on
$\Sigma^{\Box}\cup\Delta^{\Box}\cup\Gamma^{\Box}\cup\Gamma_1\cup\Gamma_2$.
We conclude that the envelope of $\Omega$ is triple sheeted. Over
$\Omega$, however, the envelope is only double sheeted,
specifically over those points in the pairwise intersections of
our three families of discs. In fact, the envelope of $\Omega$ can
be identified with the Riemann domain given as the disjoint union
of $\Sigma^{\Box}$, $\Delta^{\Box}$ and $\Gamma^{\Box}$, connected
by the sets $\Gamma_1$ and $\Gamma_2$, equipped with the natural
projection. Certainly every holomorphic function on $\Omega$
extends to this Riemann domain. By construction it is an
unbranched Riemann domain which is locally pseudoconvex. By the
equivalence of local and global pseudoconvexity of unbranched
Riemann domains (\cite[Corollary
2.2.16]{jarnicki-pflug:extension}), this Riemann domain is Stein,
and thus is the envelope of $\Omega$.

\subsection{Counterexamples in $\mathbb{C}^3$}
\label{subsection:examples in C3} We construct a domain, $\Omega$,
in $\mathbb{C}^3$ such that the envelope of $\Omega$ has one sheet
over $\Omega$, but two sheets over $\mathbb{C}^3$. This shows that
Theorem \ref{theorem:main theorem} cannot be generalized to give
information about domains in $\mathbb{C}^n$, $n>2$.\\

The key point in the example is that the intersection of one
dimensional varieties in $\mathbb{C}^3$ is not generically
preserved under
perturbation.\\

We define $\Omega$ as follows. We first build two domains, $V_1$
and $V_2$, with strictly pseudoconvex boundary points. We join
these domains with a well chosen path $\gamma$. We let $\pi_j$ be
the projection onto the $j$th coordinate of $\mathbb{C}^3$.

Begin by defining two domains, $U_1$ and $U_2$. Let $U_1$ be
\begin{align*}
U_1=
\biggl[\left\{|z|<8\right\}&\times\left\{\frac{1}{2}<|w|<1\right\}
\times\left\{|\zeta|<\frac{1}{2}\right\}\biggr]\\
& \cup \biggl[\left\{\frac{1}{2}<|z|<8\right\}
\times\left\{\frac{1}{2}<|w|<1\right\}
\times\left\{|\zeta|<1\right\}\biggr].
\end{align*}

We see that $\widetilde{U}_1$ is
\[\widetilde{U}_1=
\left\{|z|<8\right\}\times\left\{\frac{1}{2}<|w|<1\right\}\times\left\{|\zeta|<1\right\}.\]

Let $U_2$ be
\begin{align*}
U_2=
\left[\biggl\{|z|<\frac{1}{4}\right\}&\times\left\{|w|<8\right\}\times
\left\{\frac{3}{2}<|\zeta|<2\right\}\biggr] \cup\\
&\biggl[\left\{|z|<\frac{1}{4}\right\}\times\left\{\frac{3}{2}<|w|<8\right\}
\times \left\{\frac{3}{4}<|\zeta|<2\right\}\biggr].
\end{align*}

We see that $\widetilde{U}_2$ is
\[\widetilde{U}_2=\left\{|z|<\frac{1}{4}\right\}\times\left\{|w|<8\right\}
\times\left\{\frac{3}{4}<|\zeta|<2\right\}.\]

Notice that
\[U_1 \cap U_2 =\emptyset,\quad \widetilde{U}_1\cap U_2=\emptyset,
\quad \widetilde{U}_2\cap U_1=\emptyset,\] but that
$\widetilde{U}_1 \cap \widetilde{U}_2$ is non empty.

We define
\[V_1= U_1 \cap B^3(0,6)\]
and
\[V_2= U_2 \cap B^3(0,6).\]
Just as with $U_1$ and $U_2$ we have that
\[V_1 \cap V_2 =\emptyset,\quad \widetilde{V}_1\cap V_2=\emptyset,
\quad \widetilde{V}_2\cap V_1=\emptyset,\] but that
$\widetilde{V}_1 \cap \widetilde{V}_2$ is non empty.

Unlike $U_j$, however, each $V_j$ has strictly pseudoconvex
boundary points: points in $\partial V_j\cap\partial B^3(0,\,6)$.
Let $p_1$ and $p_2$ be such boundary points in $V_1$ and $V_2$,
respectively. Let $\gamma:[0,\,1]\rightarrow\mathbb{C}^3$ be a
smooth curve satisfying:
\begin{enumerate}
\item $\gamma$ runs from $p_1$ to $p_2$.
\item $\gamma$ is transversal to $\partial V_j$ at $p_j$.
\item $\gamma$ intersects $\overline{B^3(0,6)}$ only at $p_1$ and $p_2$.
\item Let $p$ be a point in the complement of $\pi_1(B^3(0,\,6))$, and define the set
\[A=\mathbb{C}^3\backslash(\{z=p\}\times\mathbb{C}^2).\]
The fundamental group of $A$ is nontrivial. We choose $\gamma$ in such a way that
\begin{enumerate}
\item $V_1\cup V_2\cup \gamma$ contains no nontrivial element of
the fundamental group, and
\item $\widetilde{V}_1\cup\widetilde{V}_2\cup\gamma$ does contain a
nontrivial element of the fundamental group.
\end{enumerate}
\end{enumerate}

As in the previous example, we find a neighbourhood, $\Gamma$, of
$\gamma$ such that $\Gamma$ is pseudoconvex, $\Gamma\cup V_1$ is
locally pseudoconvex near $p_1$, and $\Gamma\cup V_2$ is locally
pseudoconvex near $p_2$.

We define $\Omega$ as $\Omega=V_1\cup V_2 \cup \Gamma$.\\

Our choice of path ensures that $f(z,w,\zeta)=(z-p)^{1/2}$ is a
holomorphic function on $\Omega$. However, $f$ is not holomorphic
on $\widetilde{V}_1\cup\widetilde{V}_2\cup\Gamma$: as we follow
$\gamma$ from $\widetilde{V}_1$ to $\widetilde{V}_2$ we change
branches of the square root function. We conclude that $f$ is
double valued on $\widetilde{V}_1\cup\widetilde{V}_2\cup\Gamma$.
Thus the envelope of $\Omega$ is double sheeted. The two sheets
lie over $\widetilde{V}_1\cap\widetilde{V}_2$. Since this
intersection contains no points in $\Omega$, we see that the
envelope is single sheeted over $\Omega$. As in the previous
example, we are viewing the envelope as a particular unbranched
Riemann domain: the disjoint union of $\widetilde{V}_1$ and
$\widetilde{V}_1$, connected with the set $\Gamma$, and equipped
with the natural projection.

\thanks{{\em Acknowledgements}
The author would like to thank Professor Berit Stens{\o}nes for
invaluable advice in conducting this research, and Professor
Harold Boas for comments and assistance in preparation of this
article. Additionally, the referee  made valuable suggestions
which significantly improved the article. This article is part of
the author's Ph.D. dissertation.}

\nocite{pflug:pc}
\bibliographystyle{amsalpha}
\bibliography{bibliography}

\end{document}